\documentclass[moor,nonblindrev]{informs1}

\usepackage{physics}
\usepackage{mathtools}
\usepackage{tensor}

\newcommand{\half}{\tfrac{1}{2}}
\newcommand{\asymx}{\mathop{\sim}}
\newcommand{\asym}[1]{\mathrel{\asymx_{#1}}}
\newcommand{\diag}{\mathop{\rm diag} \nolimits}

\newcommand{\conv}{\mathop{\rm conv} \nolimits}

\newmuskip\pFqmuskip
\newcommand*\pFq[6][8]{%
  \begingroup 
  \pFqmuskip=#1mu\relax
  \mathcode`\,=\string"8000
  \begingroup\lccode`\~=`\
  \lowercase{\endgroup\let~}\pFqcomma
  {}_{#2}F_{#3}{\left[\genfrac..{0pt}{}{#4}{#5};#6\right]}%
  \endgroup
}
\newcommand{\pFqcomma}{\mskip\pFqmuskip}

\def\wscl{0.8}
\def\hscl{0.56}

\def\wsclx{1}
\def\hsclx{0.7}

\usepackage{natbib}
\NatBibNumeric
\def\bibfont{\small}%
\bibpunct[, ]{[}{]}{,}{n}{}{,}%

\definecolor{myblue}{RGB}{0, 20, 114}
\usepackage[colorlinks=true,breaklinks=true,bookmarks=true,urlcolor=myblue,
     citecolor=myblue,linkcolor=myblue,bookmarksopen=false,draft=false]{hyperref}

\def\EMAIL#1{\href{mailto:#1}{#1}}

\TheoremsNumberedThrough     

\EquationsNumberedThrough    

\MANUSCRIPTNO{0} 

\def\theARTICLETOP{}
\def\theLRHSecondLine{}
\def\theRRHSecondLine{}

\begin{document}


\RUNAUTHOR{Zuk and Kirszenblat}

\RUNTITLE{Non-Preemptive Priority Queue}

\TITLE{Explicit Results for the Distributions of Queue Lengths for a Non-Preemptive Two-Level Priority Queue}

\ARTICLEAUTHORS{%
\AUTHOR{Josef Zuk}
\AFF{Defence Science and Technology Group, Melbourne, Australia,
\EMAIL{josef.zuk@defence.gov.au}}
\AUTHOR{David Kirszenblat}
\AFF{Defence Science and Technology Group, Melbourne, Australia,
\EMAIL{david.kirszenblat@defence.gov.au}}
} 

\ABSTRACT{%
Explicit results are derived using simple and exact methods for the joint and marginal
queue-length distributions for the M/M/$c$ queue with two non-preemptive priority levels.
Equal service rates are assumed. Two approaches are considered.
One is based on numerically robust quadratic recurrence relations.
The other is based on a complex contour-integral representation that yields
exact closed-form analytical expressions,
not hitherto available in the literature,
that can also be evaluated numerically with very high accuracy.
}%

\MSCCLASS{Primary: 90B22; secondary: 60K25, 60J74}
\ORMSCLASS{Primary: Queues: Priority; secondary: Queues: Markovian}

\HISTORY{Date created: July 03, 2023. Last update: September 14, 2023.}

\maketitle

%

\section{Introduction}
\label{intro}
This work is concerned with the development of practical algorithms for the computation
of joint and marginal distributions of queue lengths for the M/M/$c$ queue with a
non-preemptive priority discipline.
Applications of this model are found in telecommunications \citep{NP:Cohen56},
health care \citep{NP:Hou20,NP:Taylor80},
radar \citep{NP:Orman95,NP:Orman96},
air traffic control \citep{NP:Pestalozzi64}
and numerous other areas.

The non-preemptive priority queue discipline is as stated by \citet{NP:Dressin57}:
Once a client's service has begun, it is permitted to proceed to
completion. If a server becomes empty, and there is at least one client waiting,
then a client of the highest priority present in the queue is admitted to the server.
Clients of equal priority are served on a first-come, first-served basis.
This is also known as the `head of the line' discipline.

Thus, let us consider a non-preemptive queue with K priority levels, each with a distinct
Poisson arrival rate $\lambda_k$,
\mbox{$k = 1,2,\ldots,K$}
and corresponding level traffic intensity\footnote{Consistent with \citep{NP:Gail88,NP:Kao90},
  the notation $\rho_k$ reserved for $\rho_k \equiv \lambda_k/\mu$, so that $r_k = \rho_k/N$.}
\mbox{$r_k = \lambda_k/(N\mu)$},
leading to a total traffic intensity for the aggregation of all arrivals of
\mbox{$r = \sum_{k=0}^K r_k$}.
We adopt the usual convention that smaller priority-level indices $k$
represent higher priorities. Thus, $r_1$ denotes the traffic intensity
associated with the highest priority level.
For simplicity, we have assumed a common exponential service rate $\mu$
among all priority levels.
The number of servers is denoted by
\mbox{$c = N$}.

In this work, attention is confined to the two-level problem
\mbox{$K = 2$}.
Analysis of this case is amenable to a number of analytical techniques that
do not extend easily, or at all, to the general multi-level priority problem.
Also, the two-level problem has a distinguished status.
All marginal distributions for the multi-level problem can be inferred from
the low-priority marginal pertaining to just two priority levels \citep{NP:Davis66}.
If we let $r_{\text{hi}}$ and $r_{\text{lo}}$ denote the level traffic intensities
for the high and low priority arrivals, respectively, for the two-level problem, then
the wait-conditional\footnote{See below in Section \ref{prelim}.}
marginal distribution of the queue length for priority level
\mbox{$k = 1,2,\ldots,K$}
in the multi-level problem is obtained by making the identifications
\begin{equation}
r_{\text{lo}} = r_k \;, \quad r_{\text{hi}} = \sum_{\ell = 1}^{k-1} r_\ell \;,
\end{equation}
so that the total traffic intensity in the effective (wait-conditional)
two-level problem becomes
\mbox{$r = r_{\text{sum}}$},
with
\begin{equation}
r_{\text{sum}} = r_{\text{lo}} + r_{\text{hi}} = \sum_{\ell = 1}^{k} r_\ell \;.
\end{equation}
For the actual two-level problem, we have the identifications
\mbox{$r_{\text{hi}} \equiv r_1$},
\mbox{$r_{\text{lo}} \equiv r_2$},
and we shall use both sets of notation interchangeably.
It is also convenient to introduce the parameter $\nu$ that represents the
fraction of all arrivals that are of high priority (which we abbreviate as `hifrac').
Thus
\mbox{$r_{\text{hi}} = \nu r$},
\mbox{$r_{\text{lo}} = (1-\nu) r$},
\mbox{$0 \leq \nu \leq 1$}.

Previous work on the non-preemptive priority queue has focused, almost entirely,
on calculating moments and the waiting-time distributions per priority level.
In early work, \citet{NP:Cobham54,NP:Cobham55}, followed by \citet{NP:Holley54},
were first to consider the mean waiting times and queue lengths.
Waiting-time means and second moments for general service-time distributions were
subsequently given by \citet{NP:Kesten57}.
\citet{NP:Gail88} studied the non-preemptive M/M/$c$ system for two priority levels with
different exponential service rates.
While they developed a matrix algorithm for determining various characteristics of a
generating function for this problem, explicit results were also confined to the
mean waiting times and queue lengths.

For the waiting-time problem, \citet{NP:Davis66} improved on previous work by
\citet{NP:Dressin57} to derive an explicit integral
expression for the waiting-time distribution for the non-preemptive priority queue.
He analysed the two-level problem, as the waiting time distribution for the multi-level problem
can be inferred from the two-level case.
He did not study the queue-length marginals, and they cannot be directly inferred from the
waiting-time distributions by appealing to the distributional form of
Little's law \citep{NP:Bertsimas95,NP:Keilson88}
as the no-overtaking assumption is violated.
\citet{NP:Kella85} covered the same ground as Davis for the
probability waiting function (PGF) of the waitingtime,
but using a different methodology.
The moment generating function (MGF) of the waiting time and associated moments
have also been considered in \citep{NP:Miller60}.
More recently, \citet{NP:Wagner97} has studied the waiting-time MGF for a finite-capacity,
multi-server version of the same problem as Davis.

For the queue-length distributions, \citet{NP:Miller81,NP:Miller82} uses a matrix-geometric method
for the two-level problem
that results in a complex algorithm involving multiple levels of recursion.
Little is presented about the numerical stability of this approach, and it is
known to deteriorate for traffic intensities close to unity.
\citet{NP:Kao90} and \citet{NP:Kao99} also employ matrix-geometric methods
for the two-level problem which, as they point out,
unavoidably require finite-state truncation.
The aforementioned papers deal with unequal service rates.
The matrix-geometric method \citep{NP:Neuts81} applied to queueing models
has the singular disadvantage that it necessitates truncation of the
problem to prescribed finite maximum values of queue lengths.
While powerful, it is complex and not elegant.
Thus, its use should best be avoided whenever simpler alternatives are available,
and this is manifestly the case for the present problem, as will become clear.

In earlier work, \citet{NP:Marks73} studied the two-level problem with common service rate
and derived a highly complex system of linear partial difference equations
that must be solved recursively.
The required manipulations are cumbersome and no insight into the analytic structure of the
problem is gained. However, it is most likely the first paper where actual queue-length
probabilities, rather than the PGF, were computed.
No light is shed on the numerical stability of the method.

A different approach, based on a partial PGF, is due to \citet{NP:Cohen56},
who studied the two-level problem with equal service rates;
and it is this approach that we pursue in the discussion that follows.
We take up the programme where \citet{NP:Cohen56} left off,
in devising simple and practical schemes for extracting actual probabilities from the PGF.
There is the additional benefit that this approach
can be extended to the general multi-level problem.
\citet{NP:Shortle18} have remarked that
`the determination of stationary probabilities
in a non-preemptive Markovian system is an exceedingly difficult matter, well near
impossible when the number of priorities exceeds two'.
In a separate forthcoming paper, we shall demonstrate otherwise.

The present work focuses on explicit results that are useful for practical applications.
While we do not purport to have made general theoretical advances in priority queues,
the work does serve to fill a large gap in the literature by establishing basic results
for a paradigmatic model that one would expect to have been uncovered decades ago.
We believe that it also has pedagogical value.
For the two-level non-preemptive priority queue, \citet{NP:Shortle18},
in the most recent edition of their textbook,
set up the stationary balance equations but remark that
`obtaining a reasonable solution to these stationary
equations is very difficult, ... The
most we can do comfortably is obtain expected values via two-dimensional
generating functions'.
The simplicity of the methods described herein might render
a more detailed treatment of the subject suitable for elementary texts.

\section{Non-Preemptive Priority Queue}
\label{prelim}
The no-wait probability $P_{\text{NW}}$ is the probability that a new arrival will find
at least one server idle. It is clearly independent of the queue discipline, and
is given by \citep{NP:Davis66}
\begin{equation}
\frac{1}{1 - P_{\text{NW}}} = 1 + (1-r)\frac{N!}{(Nr)^N}{\cdot}
     \sum_{k=0}^{N-1} \frac{(Nr)^k}{k!} \;.
\label{PNW}
\end{equation}
Let
\mbox{$P(n,m)$}
denote the steady-state probability that there are $n$ low-priority clients in the queue
(rather than in the system) and $m$ high-priority clients in the queue.
We have the decomposition \citep[\mbox{\it cf.}][]{NP:Davis66}
\begin{equation}
P(n,m) = P_{\text{NW}}{\cdot}\delta_{n0}\delta_{m0} +
     (1 - P_{\text{NW}}){\cdot}f(n,m) \;,
\end{equation}
where $f(n,m)$ represents the wait-conditional joint PMF, {\it i.e.}\
the probability that there are $n$ low-priority clients and $m$ high-priority
clients in the queue, given that all servers are busy.
The wait-conditional distribution does not explicitly depend on the number
of servers $N$.
There is only an indirect dependence on $N$ through the total traffic intensity $r$.

Our starting point is the paper of \citet{NP:Cohen56},
which introduced a partial PGF for the problem that summed only
over the low-priority argument:
\begin{equation}
G_m(p) \equiv \sum_{n=0}^\infty p^n P(n,m) \;.
\end{equation}
This turns out to be a very convenient strategy, especially given the fact that
the wait-conditional high-priority marginal is a simple geometric distribution.
Only the low-priority marginal is non-trivial.
We introduce a wait-conditional version $g_m(p)$ of this PGF such that
\begin{equation}
g_m(p) \equiv \sum_{n=0}^\infty p^n f(n,m) \;.
\end{equation}
It follows that
\begin{equation}
G_m(p) = P_{\text{NW}}{\cdot}\delta_{m0} +
     (1 - P_{\text{NW}}){\cdot}g_m(p) \;.
\end{equation}

Cohen's result \citep{NP:Cohen56} for the wait-conditional PGF
for the two-level non-preemptive priority queue with equal service rates
is\footnote{The quantities $\lambda_{1,2}(p)$ should not be confused with the arrival rates introduced earlier.
We are adhering to Cohen's original, but less than ideal, notation.}
\begin{equation}
g_m(p) = \frac{(1-r)(1-p)}{1 - p\lambda_2(p)}{\cdot}\lambda_1^m(p) \;,
\label{gCohen}
\end{equation}
where
\mbox{$\lambda_{1,2}(p)$}
are defined as follows:
Let us introduce
\mbox{$\lambda(p) = \lambda_\pm(p)$}
as the two solutions of the quadratic equation
\begin{equation}
\lambda^2 - (1 + r - r_2 p)\lambda + r_1 = 0 \;,
\end{equation}
such that
\begin{equation}
\lambda_\pm(p) = [b(p) \pm \sqrt{b^2(p) - 4r_1}]/2 \;, \quad
     b(p) \equiv 1 + r - r_2 p \;.
\end{equation}
Then in (\ref{gCohen}), we have
\mbox{$\lambda_1(p) = \lambda_-(p)$},
\mbox{$\lambda_2(p) = \lambda_+(p)$},
and it is useful to note that
\begin{equation}
\lambda_+(p)  + \lambda_-(p) = b(p) \;, \quad
     \lambda_+(p)\cdot\lambda_-(p) = r_1 \;.
\end{equation}

Another way to express the
PGF for the wait-conditional distribution is
\begin{equation}
g_m(p) = g_{\text{lo}}(p){\cdot}[1 - \lambda_1(p)]\lambda_1^m(p) \;, \quad
     g_{\text{lo}}(p) = \frac{1 - r}{\lambda_2(p) - r} \;.
\label{gCohen1}
\end{equation}
That
\mbox{$g_{\text{lo}}(p)$}
represents the wait conditional PGF for the low-priority marginal is clear from
observing that
\begin{equation}
\sum_{m = 0}^\infty g_m(p) = g_{\text{lo}}(p) \;.
\end{equation}
On the other hand, it follows directly from (\ref{gCohen1})
that the wait-conditional PMF for the high priority marginal is given by
\begin{equation}
f_{\text{hi}}(m) = g_m(1) = (1 - r_1)r_1^m \;.
\end{equation}
Consequently, the only marginal distribution of interest in the present study is that for the
low-priority level.

By construction, the wait-conditional joint PMF $f(n,m)$ is recovered from the PGF $g_m(p)$
according to
\begin{equation}
f(n,m) = \frac{1}{n!}{\cdot}\left.\frac{d^n}{dp^n} g_m(p)\right|_{p=0} \;.
\label{PGFDeriv}
\end{equation}
The multiple derivative is prohibitively cumbersome to directly perform analytically.
Thus, we proceed to present two alternative strategies that render the problem tractable.

\section{Quadratic Recurrence}
The first method constructs a recurrence relation based on the fact that the functions
$\lambda_\pm(p)$ solve a quadratic equation.
We begin by considering the low-priority marginal, whose PGF can be expressed as
\begin{equation}
g_{\text{lo}}(p) = \frac{1 - r}{\lambda_2(p) - r} \;.
\end{equation}
Since $\lambda_2(p)$ satisfies a quadratic equation, then so does $g_{\text{lo}}(p)$.
Let us set
\begin{equation}
u \equiv r_2p \;, \quad g(u) \equiv \frac{1}{\lambda_2 - r} = \sum_{k= 0}^\infty \frac{g_k}{k!} u^k
     \quad \Rightarrow
     \quad g_k = \left.\frac{d^k g(u)}{du^k}\right|_{u=0} \;.
\end{equation}
Then we obtain
\begin{equation}
(ru - r_2)g^2 + (u-1+r)g + 1 = 0 \;.
\end{equation}
We now differentiate this equation $n$ times with respect to $u$, and use the identities
\begin{align}
\begin{aligned}
\frac{1}{n!}{\cdot}\left.\frac{d^n}{du^n} (ug)\right|_{u=0}   &= \frac{g_{n-1}}{(n-1)!} \;, \\
\frac{1}{n!}{\cdot}\left.\frac{d^n}{du^n} (g^2)\right|_{u=0}  &= \sum_{k=0}^n \frac{g_{k}}{k!}{\cdot}\frac{g_{n-k}}{(n-k)!} \;, \\
\frac{1}{n!}{\cdot}\left.\frac{d^n}{du^n} (ug^2)\right|_{u=0} &= \sum_{k=0}^{n-1} \frac{g_{k}}{k!}{\cdot}\frac{g_{n-k-1}}{(n-k-1)!} \;.
\end{aligned}
\end{align}
For the quantities
\mbox{$f_k \equiv g_k/k!$},
this leads to the non-linear recurrence relations
\begin{equation}
f_n = \frac{1 + rf_0}{1 - r + 2r_2f_0}{\cdot}f_{n-1} +
     \frac{1}{1 - r + 2r_2f_0}\sum_{k=1}^{n-1} f_k{\cdot}(rf_{n-k-1} - r_2f_{n-k}) \;,
\end{equation}
for
\mbox{$n = 1,2,\ldots$},
with
\begin{equation}
f_0 = \frac{1}{2r_2}\left[\sqrt{(1-r)^2 + 4r_2} - (1-r)\right]
    = \frac{2}{1 - r + \sqrt{(1-r)^2 + 4r_2}}
    > 0\;.
\end{equation}
The expression for $f_0$ follows from
\mbox{$f_0^{-1} = g_0^{-1} = \lambda_2(0) - r$}.
We observe that
\mbox{$f_{\text{lo}}(n) = (1-r)r_2^n f_n$}.

Efficient vectorized implementations in {\sc Matlab} are possible.
Practical implementation proceeds as follows:
Let us introduce an arbitrary scale factor $\Lambda$, define
\begin{equation}
D \equiv 1 - r + 2r_2f_0 = \sqrt{(1-r)^2 + 4r_2} \;,
     \quad c_1 \equiv r_2/\Lambda \;,
     \quad c_2 \equiv \Lambda/D \;,
\end{equation}
and scale according to
\mbox{$\tilde{f}_n \equiv \Lambda^n f_n = (r_2/c_1)^n f_n$}.
Then we solve the recurrence
\begin{equation}
\tilde{f}_n = c_2{\cdot}\biggl(\tilde{f}_{n-1} + \sum_{k=0}^{n-1} \tilde{f}_k \Delta_{(n-1)-k}\biggr)
\label{SclRecur}
\end{equation}
and recover the marginal as
\mbox{$f_{\text{lo}}(n) = (1-r)c_1^n\tilde{f}_n$}.
At each step, we set
\begin{equation}
\Delta_k \equiv r\tilde{f}_k - c_1\tilde{f}_{k+1} \;,
\end{equation}
for
\mbox{$k = 1,2,\ldots,n-1$},
subject to the initialization
\mbox{$\tilde{f}_n \leftarrow 0$}
within the scope of evaluating $\Delta_{n-1}$.
We find that good numerical performance is achieved with
\mbox{$\Lambda = r_2$},
so that
\mbox{$c_1 = 1$}.

Analogous treatment of the joint PMF is only marginally more complex.
Based on the quadratic
\begin{equation}
\lambda_\pm^2 + (u - 1- r)\lambda_\pm +r_1 = 0\;,
\end{equation}
we solve for the Taylor-series coefficients $\lambda^{(k)}_\pm$ in
\begin{equation}
\lambda_\pm = \sum_{k=0}^\infty \lambda^{(k)}_\pm p^k = \sum_{k=0}^\infty \Lambda^{-k} f^\pm_k u^k\;,
\end{equation}
for some arbitrary scale factor $\Lambda$,
using the non-linear recurrence
\begin{equation}
f^\pm_n = \mp\frac{1}{\sqrt{(1-r)^2 + 4r_2}}\left(\Lambda f^\pm_{n-1} +
     \sum_{k=1}^{n-1} f^\pm_k{\cdot}f^\pm_{n-k}\right) \;,
\end{equation}
\mbox{$n = 1,2,\ldots$},
where
\begin{equation}
f^\pm_0 = \half\left(1 + r \pm\sqrt{(1-r)^2 + 4r_2}\right) \;.
\end{equation}
The $\lambda$-coefficients are recovered according to
\mbox{$\lambda^{(k)}_\pm = (r_2/\Lambda)^k f^\pm_k$}.
As with the marginal, the choice
\mbox{$\Lambda = r_2$}
results in good numerical performance.
All that remains to be done is to use the standard recursion for multiplication
of power series as dictated by (\ref{gCohen1}).
The simplest way to proceed is via repeated convolutions:
\begin{align}
\begin{aligned}
\phi_0 &= (1-r){\cdot}\conv\left(\frac{1}{\lambda_2 - r}, 1-\lambda_1\right) \;, \\
\phi_k &= \conv(\phi_{k-1}, \lambda_1) \;,
\end{aligned}
\end{align}
for
\mbox{$k = 1,2,\ldots,m$}.
Then
\mbox{$f(n,m) = \phi_m(n)$}.
The conv function is defined like
the {\sc Matlab} function of the same name:
Suppose that
\mbox{$C(u) = A(u)B(u)$},
with
\begin{equation}
A(u) = \sum_{n=0}^{n_1} a(n)u^n \;, \quad
    B(u) = \sum_{n=0}^{n_2} b(n)u^n \;, \quad
    C(u) = \sum_{n=0}^{n_1+n_2} c(n)u^n \;.
\end{equation}
Then
\mbox{$c = \conv(a,b)$},
where
\begin{equation}
c(n) = \conv(a,b)(n) \equiv \sum_{k=0}^{n} a(k)b(n-k) \;,
\end{equation}
for
\mbox{$n = 0,1,2,\dots,n_1+n_2$}.

The foregoing recurrence relations constitute a significant improvement over the strategy
implemented in \citep{NP:Cidon90},
and are vastly simpler than those arising from the matrix-geometric method as considered
in \citep{NP:Kao90,NP:Kao99,NP:Miller81,NP:Miller82}.
While the quadratic recurrence method exhibits excellent numerical behaviour,
it gives little insight into the analytical structure of the distributions.
This deficiency is addressed in the next section.

\section{Complex Contour Integral}
\label{CCIntegral}
\begin{figure}
\FIGURE
{\includegraphics[width=\wscl\linewidth, height=\hscl\linewidth]{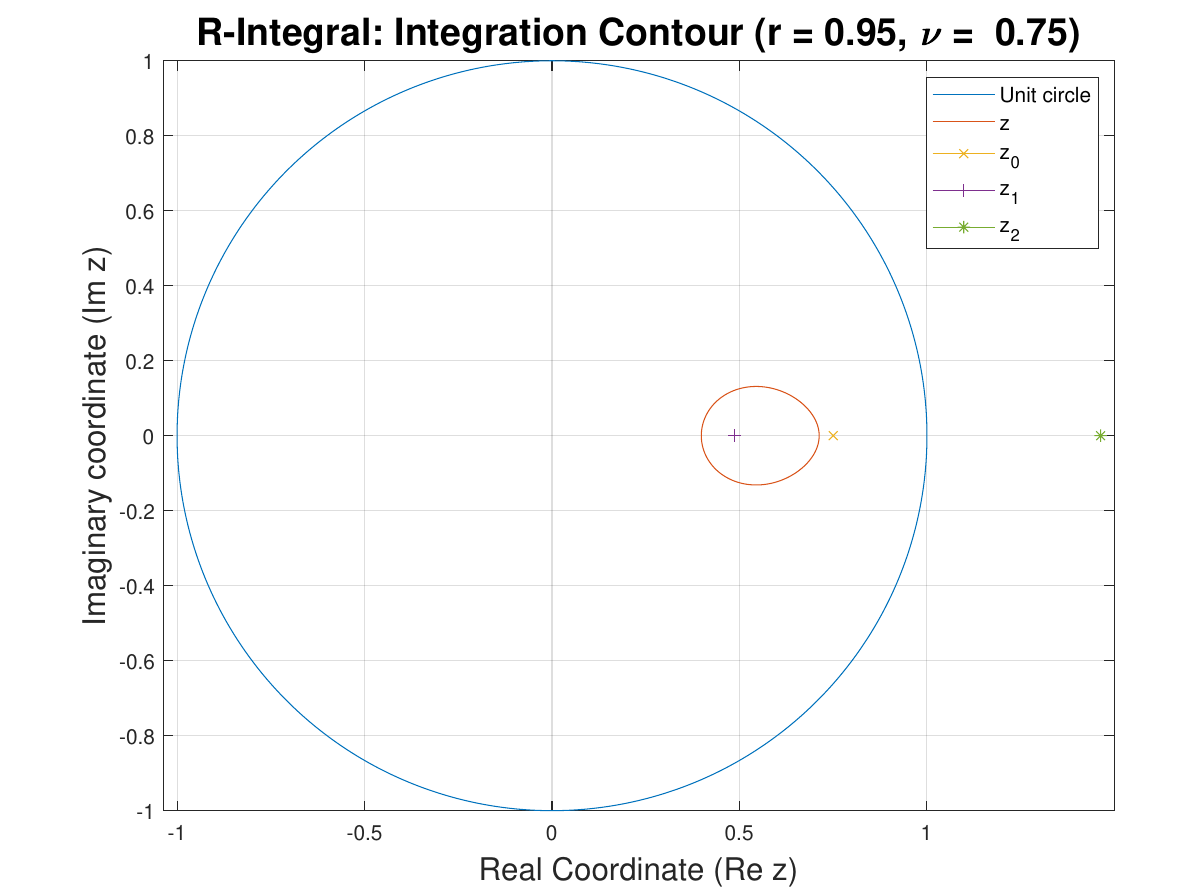}}
{\hphantom{x}\label{ZContour}}
{The $z$-contour that results from taking the $p$-contour to be the unit circle centred
     on the origin, plotted for the case of
     total traffic intensity $r = 0.95$ and fraction of high-priority arrivals $\nu = 0.75$.
     Also displayed are the locations of R-integral poles $z_0,z_1, z_2$.}
\end{figure}

Another strategy in dealing with (\ref{PGFDeriv}) is to represent it
in terms of a complex contour integral
in accordance with Cauchy's integral theorem.
This yields
\begin{equation}
f(n,m) = (1-r)\oint_{\mathcal{C}}\frac{dp}{2\pi i}\, \frac{(1-p)\lambda_1^m}{p^{n+1}(1 - p\lambda_2)} \;,
\label{WCPMF}
\end{equation}
where $\mathcal{C}$ is an anti-clockwise circle centred about the origin with radius less than $1/r$.
It follows directly that the
low-priority marginal PMF, defined by
\begin{equation}
f_{\text{lo}}(n) \equiv \sum_{m= 0}^\infty f(n,m) \;,
\end{equation}
is represented as a complex contour integral by
\begin{equation}
f_{\text{lo}}(n) = (1-r)\oint_{\mathcal{C}}\frac{dp}{2\pi i}\, \frac{1-p}{p^{n+1}}{\cdot}
     \frac{1}{(1 - p\lambda_2)(1 - \lambda_1)} \;.
\end{equation}

The conventional approach in dealing with such contour integrals, mirroring the approach
adopted previously for the waiting-time distribution \citep{NP:Davis66}, would be to deform
the contour by expanding it to the circle at infinity while avoiding a cut of finite extent
on the real axis that is generated by the square-root component of $\lambda_\pm(p)$,
and a possible simple pole that also lies on the real axis.
The circle at infinity yields a vanishing contribution, which leaves a (potential) pole term
and a real-valued integral along the cut.
We shall explore this approach separately in a forthcoming paper, where we shall show that
it leads to integral expressions that are amenable to efficient quadrature algorithms, and can
also be evaluated analytically in terms of a generalized form of the associated Legendre functions.
In the present work, we pursue a different method based on a change of integration variable.

Let
\mbox{$\lambda = z_1,z_2$}
be the roots of the polynomial equation
\mbox{$\lambda^2 - (1+r)\lambda + r_1 = 0$},
so that we have
\mbox{$z_1 + z_2 = 1 + r$},
\mbox{$z_1z_2 = r_1 = \nu r$}.
Then, the inversion of
\mbox{$z = \lambda_\pm(p)$}
yields
\begin{equation}
p = -(z - z_1)(z - z_2)/(r_2z) \;.
\end{equation}
Thus,
\begin{equation}
dp = -\frac{1}{r_2}\left(1 - \frac{z_1z_2}{z^2}\right) dz \;,
\end{equation}
in which case
\begin{equation}
\frac{dp}{p^{n+1}} = (-r_2)^n\frac{(z^2 - z_1z_2)z^{n-1}}
     {\left[(z-z_1)(z-z_2)\right]^{n+1}}{\cdot}dz \;.
\end{equation}
We make the change of integration variable
\mbox{$p\mapsto z: z = \lambda_1(p)$},
in which case
\mbox{$\lambda_2(p) = r/z$},
and we make the identifications
\begin{equation}
z_0 = r_1/r \;, \quad z_1 = \lambda_-(p = 0) \;, \quad z_2 = \lambda_+(p = 0) \;,
\end{equation}
or, equivalently,
\begin{equation}
z_0 = \nu \;, \quad
z_1 = \half\left[1 + r - \sqrt{(1 + r)^2 - 4\nu r}\right] \;, \quad
z_2 = \half\left[1 + r + \sqrt{(1 + r)^2 - 4\nu r}\right] \;.
\label{z012}
\end{equation}
Then, we obtain
\begin{equation}
\frac{1 - p}{1 - p\lambda_2} = \frac{z}{r}{\cdot}\frac{z - 1}{z - z_0} \;.
\end{equation}
It follows that the joint PMF is given by
\begin{equation}
    f(n,m) = \frac{(1-r)(-r_2)^n}{r}\oint_{\mathcal{C}'}\frac{dz}{2\pi i}\, \frac{z^{m+n}}{z - z_0}{\cdot}
    \frac{(z - 1)(z^2 - z_1z_2)}{[(z - z_1)(z - z_2)]^{n+1}} \;,
\end{equation}
where $\mathcal{C}'$ is a closed anti-clockwise contour
that encloses the pole at
\mbox{$z = z_1$}.
but with the poles at
\mbox{$z = z_0, z_2$}
in the exterior.
For the low-priority marginal PMF, we have
\begin{equation}
    f_{\text{lo}}(n) = -\frac{(1-r)(-r_2)^n}{r}\oint_{\mathcal{C}'}\frac{dz}{2\pi i}\, \frac{z^n}{z - z_0}{\cdot}
    \frac{(z^2 - z_1z_2)}{[(z - z_1)(z - z_2)]^{n+1}} \;.
\end{equation}

\section{R-Integrals}
In order to evaluate the integral representations for the joint and marginal PMFs, derived the
foregoing section,
we introduce a collection of complex contour integrals, to which we shall refer as the R-integrals,
according to the definition
\begin{equation}
R^m_n \equiv \oint_{\mathcal{C}'} \frac{dz}{2\pi i}\, \frac{1}{z - z_0}{\cdot}
     \frac{z^m}{[(z - z_1)(z - z_2)]^n} \;,
\label{RIntg}
\end{equation}
for
\mbox{$m,n = 0,1,2,\ldots$},
where $\mathcal{C}'$ is a closed anti-clockwise contour
that encloses the pole at
\mbox{$z = z_1$}.
but with the poles at
\mbox{$z = z_0, z_2$}
in the exterior.
An immediate consequence of this definition is the (backwards) recurrence relation
\begin{equation}
R^m_{n-1} = R^{m+2}_{n} - (z_1 + z_2)R^{m+1}_{n} + z_1z_2 R^{m}_{n} \;.
\label{Recur1}
\end{equation}
One may also note the scaling behaviour
\begin{equation}
R^m_n(z_0,z_1,z_2) = z_0^{m - 2n}R^m_n(1,z_1/z_0,z_2/z_0) \;,
\end{equation}
or, more generally,
\begin{equation}
R^m_n(z_0,z_1,z_2) = \zeta^{m-2n}R^m_n(z_0/\zeta,z_1/\zeta,z_2/\zeta) \;,
\end{equation}
for any
\mbox{$\zeta > 0$}.

In the present application to the priority queue, the parameters
\mbox{$z_0,z_1,z_2$}
are given by (\refeq{z012}).
In terms of the R-integrals, the joint PMF is given by
\begin{equation}
f(n,m) = \frac{(1-r)(-r_2)^n}{r}\left(R^{m+n+3}_{n+1} - R^{m+n+2}_{n+1}
    - z_1z_2R^{m+n+1}_{n+1} + z_1z_2R^{m+n}_{n+1}\right)\;.
\end{equation}
If we introduce the difference functions
\mbox{$\Delta R^m_n \equiv R^{m+1}_n - R^m_n$},
then we can write
\begin{equation}
f(n,m) = \frac{(1-r)(-r_2)^n}{r}\left(\Delta R^{m+n+2}_{n+1}
     - z_1z_2\Delta R^{m+n}_{n+1}\right)\;.
\label{fjnt}
\end{equation}
Likewise, in terms of the R-integrals, we have
for the low-priority marginal PMF,
\begin{equation}
f_{\text{lo}}(n) = -\frac{(1-r)(-r_2)^n}{r}\left(R^{n+2}_{n+1} - z_1z_2R^{n}_{n+1}\right) \;,
\label{flo}
\end{equation}
for
\mbox{$n = 0,1,2,\ldots$}.
For the exclusively-low distribution, defined by
\mbox{$f_{\text{xlo}}(n) \equiv f(n,0)$},
we can write
\begin{equation}
f_{\text{xlo}}(n) = \frac{(1-r)(-r_2)^n}{r}\left(\Delta R^{n+2}_{n+1} - z_1z_2\Delta R^{n}_{n+1}\right) \;.
\label{fxlo}
\end{equation}
It gives the probability of finding $n$ low-priority clients in the queue and no high-priority clients.
It has a form that is similar to the low-priority marginal $f_{\text{lo}}(n)$, and we will show later
that the two are, in fact, closely related. This relationship will provide a useful diagnostic test
of the numerical performance of the R-integral computation.

We have succeeded in recasting the problem into one that involves complex contour integration
over a collection of totally meromorphic functions.
In Figure~\ref{ZContour}, we plot
the $z$-contour $\mathcal{C}'$ that results from taking the $p$-contour $\mathcal{C}$
to be the unit circle centred
on the origin, plotted for the case of
total traffic intensity $r = 0.95$ and fraction of high-priority arrivals $\nu = 0.75$.
Also displayed are the locations of R-integral poles $z_0,z_1, z_2$.

\subsection{Recurrence}
If we cast the recurrence relation (\ref{Recur1}) as
\begin{equation}
R^m_{n} = R^{m-2}_{n-1} - z_1z_2 R^{m-2}_{n} + (z_1 + z_2)R^{m-1}_{n} \;,
\label{Recur2}
\end{equation}
for
\mbox{$m = 2,3,\ldots$},
\mbox{$n = 1,2,\ldots$},
then it may, in principle, be solved recursively for the $R^m_n$
starting from the seed values
\begin{align}
\begin{aligned}
R^m_{0\hphantom{+1}} &= 0 \;, \\
R^0_{n+1}            &= \frac{(-1)^n}{\left[(z_1-z_0)(z_1-z_2)\right]^{n+1}}
                        {\cdot}p_n\left(\frac{z_1 - z_0}{z_1 - z_2}\right) \;, \\
R^1_{n+1}            &= \frac{(-1)^n}{(z_1 - z_2)^{2n+1}}\binom{2n}{n} + z_0 R^0_{n+1} \;,
\end{aligned}
\end{align}
with the polynomials $p_n(x)$ defined by
\begin{equation}
p_n(x) \equiv \sum_{k = 0}^n \binom{k+n}{k} x^k \;.
\end{equation}
Unfortunately, this recursion scheme is numerically unstable,
especially for small $\nu$.

\subsection{Series Representation}
Applying Cauchy's theorem  to (\ref{RIntg}), followed by an invocation of Leibniz's formula,
we obtain
\begin{align}
\begin{aligned}
R^m_{n+1} &= \frac{1}{n!}{\cdot}\frac{d^n}{dz_1^n}\left[
     \frac{1}{(z_1-z_2)^{n+1}}{\cdot}\frac{z_1^m}{z_1-z_0}\right] \\
&= \sum_{k=0}^n \frac{1}{(n-k)!}\frac{d^{n-k}}{dz_1^{n-k}}
     \left[\frac{1}{(z_1-z_2)^{n+1}}\right]{\cdot}
     \frac{1}{k!}\frac{d^k}{dz_1^k}
     \left[\frac{z_1^m}{z_1-z_0}\right] \;.
\end{aligned}
\end{align}
The first differentiation is trivial to perform, yielding
\begin{equation}
R^m_{n+1} = (-1)^n\sum_{k=0}^n\binom{2n-k}{n}\frac{z_0^{m-k-1}}{(z_2-z_1)^{2n+1-k}}
     S^m_k(z_1/z_0) \;,
\end{equation}
where
\begin{equation}
S_k^m(x) \equiv \frac{1}{k!}\frac{d^k}{dx^k}\left(\frac{x^m}{1-x}\right) \;.
\label{SFn}
\end{equation}
The functions $S_k^m(x)$ satisfy the relationship
\begin{equation}
S_k^{m+1}(x) - S_k^m(x) = -\binom{m}{k}x^{m-k} \;.
\end{equation}
It is convenient to introduce polynomials
\begin{equation}
P_k^m(x) \equiv (1-x)^{k+1}S_k^m(x) \;,
\label{PS}
\end{equation}
so that
\mbox{$P_k^0(x) = 1$},
\mbox{$P_0^m(x) = x^m$}.
Then we can write
\begin{equation}
R^m_{n+1} = \frac{(-1)^n}{(z_2-z_1)^{2n+2}}\sum_{k=0}^n\binom{2n-k}{n}
     \left(\frac{z_2-z_1}{1 - z_1/z_0}\right)^{k+1}z_0^{m-k-1}
     P^m_k(z_1/z_0) \;.
\end{equation}
Combining (\ref{SFn}) and ({\ref{PS}),
we can establish that, for
\mbox{$m > k$},
\begin{align}
\begin{aligned}
P_k^m(x) &= (1-x)^{k+1}\sum_{\ell=0}^k \frac{1}{\ell !}\frac{d^\ell x^m}{dx^\ell}{\cdot}
     \frac{1}{(k-\ell)!}\frac{d^{k-\ell}}{dx^{k-\ell}}\left(\frac{1}{1-x}\right) \\
&= \sum_{\ell=0}^k D_\ell^m(x) \;,
\end{aligned}
\end{align}
where
\begin{equation}
D_\ell^m(x) \equiv \binom{m}{\ell}x^{m-\ell}(1 - x)^\ell \;.
\label{DPoly}
\end{equation}
Equation (\ref{DPoly}) represents a cumulative sum, each term of which can be computed recursively.
For example, when $x$ is bounded away for zero,
\begin{equation}
D^m_\ell(x) = \left(\frac{m+1}{\ell} - 1\right){\cdot}\left(\frac{1}{x} - 1\right)D^m_{\ell-1}(x) \;,
\end{equation}
for
\mbox{$\ell = 1,2,\ldots$},
with
\mbox{$D^m_0(x) = x^m$}.
A similar recursion holds for small $x$, computed backwards from
\mbox{$D^m_m(x) = (1 - x)^m$}.

An explicit representation of the polynomials $P_k^m(x)$ is given by
\begin{equation}
P_k^m(x) = 1 - (1-x)^{k+1}\sum_{\ell=0}^{m-k-1}\binom{k+\ell}{\ell} x^\ell \;.
\end{equation}
It may be observed that
\mbox{$P_k^m(x) = 1$}
whenever
\mbox{$m \leq k$},
and that
\mbox{$P_k^m(x) \geq 0$}
for all
\mbox{$0 \leq x \leq 1$}.
These polynomials also satisfy the recurrence relation
\begin{equation}
P_k^{m+1}(x) = P_k^m(x) + \frac{m}{k}(1-x)\left[P_{k-1}^m(x) - P_{k-1}^{m-1}(x)\right] \,
\end{equation}
for
\mbox{$k,m = 1,2,\ldots$},
subject to
\begin{equation}
P_0^m(x) = x^m \;, \quad P_k^0(x) = 1 \;, \quad P_k^1(x) = 1 - (1-x)\delta_{k0} \;.
\end{equation}

\subsection{Evaluation}
In order to achieve good numerical behaviour as
\mbox{$\nu\to 1$},
it is convenient to work with the scaled integrals
\mbox{$\hat{R}^m_{n+1} \equiv (-r_{\text{lo}})^n R^m_{n+1}$},
for which we have the well-behaved series representation
\begin{equation}
\begin{split}
\hat{R}^m_{n+1} = &\frac{1}{(z_2-z_1)(1 - z_1/z_0)}{\cdot}
     \left(\frac{r_{\text{lo}}}{(z_2-z_1)^2}\right)^n \\
     &{}\times\sum_{k=0}^n\binom{2n-k}{n}
     \left(\frac{z_2-z_1}{1 -z_1/z_0}\right)^k z_0^{m-k-1}P_k^m(z_1/z_0) \;.
\label{RSeries}
\end{split}
\end{equation}
Thus, we consider the computation of the vectors
\begin{equation}
\hat{\mathbf{R}}^{(m)} \equiv [\hat{R}^m_1, \hat{R}^m_2,\ldots, \hat{R}^m_{N+1}]^{\sf T} \;.
\end{equation}
To assist with this, we define the constant
\begin{equation}
\kappa \equiv \frac{1}{(z_2-z_1)(1 - z_1/z_0)} \;,
\end{equation}
the diagonal matrices
\begin{align}
\begin{aligned}
A &\equiv \diag[a^0,a^1,\ldots,a^N] \;, &  a &\equiv r_{\text{lo}}/(z_2-z_1)^2 \;, \\
B &\equiv \diag[b^0,b^1,\ldots,b^N] \;, &  b &\equiv (z_2-z_1)/(1 -z_1/z_0) \;,
\end{aligned}
\end{align}
and the combinatorial matrix
\begin{equation}
C_{nk} \equiv \binom{2n-k}{n}
\end{equation}
provided $k \leq n$ and is zero otherwise.
We also introduce the polynomial vectors
\begin{equation}
\mathbf{P}^{(m)} \equiv [P^{(m)}_0, P^{(m)}_1,\ldots, P^{(m)}_N]^{\sf T} \;, \quad
     P^{(m)}_k \equiv z_0^{m-k-1}P_k^m(z_1/z_0) \;.
\end{equation}
Then, we can write (\ref{RSeries}) as
\begin{align}
\begin{aligned}
\hat{\mathbf{R}}^{(m)} &= \kappa{\cdot}ACB\mathbf{P}^{(m)} \\
&= \kappa{\cdot}(ACA^{-1}){\cdot}AB{\cdot}\mathbf{P}^{(m)} \;.
\end{aligned}
\end{align}
At this point, we note that the product $AB$ is the diagonal matrix of increasing powers
\begin{equation}
AB = \diag[\gamma^0,\gamma^1,\ldots,\gamma^N] \;, \quad
     \gamma \equiv r_{\text{lo}}/[(z_2-z_1)(1-z_1/z_0)] \;,
\end{equation}
and that
\mbox{$(ACA^{-1})_{nk} = a^{n-k}C_{nk}$},
which is easily computed by observing the cumulative product form
\begin{equation}
a^\ell\binom{n+\ell}{n} = \prod_{j=1}^\ell\left(1 + n/j\right) a\;.
\end{equation}
If we combine the column vectors $\hat{\mathbf{R}}^{(m)}$ and $\mathbf{P}^{(m)}$
into respective matrices, so that
\begin{align}
\begin{aligned}
\hat{\mathbf{R}} &\equiv [\hat{\mathbf{R}}^{(0)}, \hat{\mathbf{R}}^{(1)}, \ldots, \hat{\mathbf{R}}^{(M)}] \;, \\
\mathbf{P}       &\equiv [\mathbf{P}^{(0)}, \mathbf{P}^{(1)}, \ldots, \mathbf{P}^{(M)}] \;,
\end{aligned}
\end{align}
then we obtain the  matrix equation
\begin{equation}
\hat{\mathbf{R}} = \kappa{\cdot}(ACA^{-1}){\cdot}AB{\cdot}\mathbf{P} \;.
\end{equation}
In Figure~\ref{LoMargPMFMulti}, we plot
the queue-length PMF for the low-priority arrivals,
as the negative base-10 logarithm, for
total traffic intensity $r = 0.99$ and a range of hifrac values $\nu$.
Overlaid, are the asymptotic curves in the large queue-length limit.
This is given by
\begin{equation}
f_{\text{lo}}(n) \asym{n\to\infty} \sqrt{\frac{1-r}{\pi r}}{\cdot}\frac{r^n}{\sqrt{n}} \;,
\end{equation}
when
\mbox{$r_\text{hi} = r^2$}
(or equivalently
\mbox{$\nu = r$}).
Otherwise, the low-priority marginal PMF can be decomposed into two components
according to
\begin{equation}
f_{\text{lo}}(n) = f_{\text{pol}}(n){\cdot}\Theta(r^2 - r_\text{hi}) + f_{\text{cut}}(n) \;,
\end{equation}
where $\Theta(x)$ denotes the Heaviside function such that
\mbox{$\Theta(x) = 1$}
for
\mbox{$x \geq 0$}
and vanishes otherwise.
The large-$n$ behaviour of these components is given by
\begin{align}
\begin{aligned}
f_{\text{pol}}(n) &\asym{n\to\infty} \left[1 - \frac{r(1-r)}{r_\text{lo}}\right](1-r)r^{n-1} \;, \\
f_{\text{cut}}(n) &\asym{n\to\infty} \frac{(\sqrt{r_\text{hi}}/r_\text{lo})^{1/2}}{2\sqrt{\pi}r}
     {\cdot}\frac{1-r}{(\chi-1/r)\chi^{n-1/2}n^{3/2}} \;,
\end{aligned}
\end{align}
where
\mbox{$\chi \equiv 1 + (1 - \sqrt{r_\text{hi}})^2/r_\text{lo} > 1/r$}.
The derivation of these results, which will be presented in a forthcoming paper, follows directly from
the pole/cut integral representation of the distribution, mentioned in Section~\ref{CCIntegral}.
The computed points, represented by the coloured dots, are interpolated by black curves.
The asymptotic curves are indicated by a coloured dashed line-style.
Thus, when the interpolation between the data points becomes coloured, this indicates that the
agreement between the computation and asymptotic limit is within the linewidth of the graph.
In Figure~\ref{LoMargPMF}, we plot
the queue-length PMF for the low-priority arrivals,
as the negative base-10 logarithm, for the case of
total traffic intensity $r = 0.99$ and fraction of high-priority arrivals $\nu = 0.95$,
where asymptotic behaviour is slow to set in.
We see that the computation remains robust up to a queue length of at least $n = 1000$
which lies deep in the asymptotic region.
In Figure~\ref{PMFMap}, we plot a
two-dimensional map of the joint probability distribution $f(n,m)$ of the queue lengths,
for total traffic intensity $r = 0.75$ and fraction of high-priority arrivals $\nu = 0.9$.
A logarithmic scaling has been applied, such that
\mbox{$f(n,m) \leftarrow \max\{0, 1 + \log_{10}(f(n,m)/f_{\text{max}})/20\}$},
where
\mbox{$f_{\text{max}} \equiv \max\{f(n,m)\}$}.

\subsection{Limiting Cases}
The
\mbox{$\nu\to 0$}
limiting behaviour of the R-integrals is given by
\begin{equation}
R_n^m \asym{\nu\to 0^+} \left\{
\begin{array}{cll}
0                                               & \quad\text{for}\quad & m > n \\
(-1)^{n-1}\left[1 - r^{n-1}/(1+r)^{2n-1}\right] & \quad\text{for}\quad & m = n \\
(-1)^{n-1}/\nu^{n-m}                            & \quad\text{for}\quad & m < n
\end{array}
\right. \;.
\label{HifracZeroLimit}
\end{equation}
At the opposite extreme, for
\mbox{$\nu = 1$},
we have
\mbox{$z_0 = 1,\; z_1 = r,\; z_2 = 1$},
in which case
\begin{equation}
\frac{z_2 - z_1}{1 - z_1/z_0} = 1 \;.
\end{equation}
It follows that
\begin{equation}
R_{n+1}^m \asym{\nu\to 1^-} \frac{(-1)^n}{(1-r)^{2n+2}}\sum_{k=0}^n \binom{2n-k}{n}
     P_k^m(r) \;.
\end{equation}
Equation (\ref{HifracZeroLimit}) shows that the R-integrals become singular for small $\nu$
when $m < n$. This is one reason for the numerical instability of the recurrence relations (\ref{Recur2}),
given that the seed values always reside in this region.

\begin{figure}
\FIGURE
{\includegraphics[width=\wscl\linewidth, height=\hscl\linewidth]{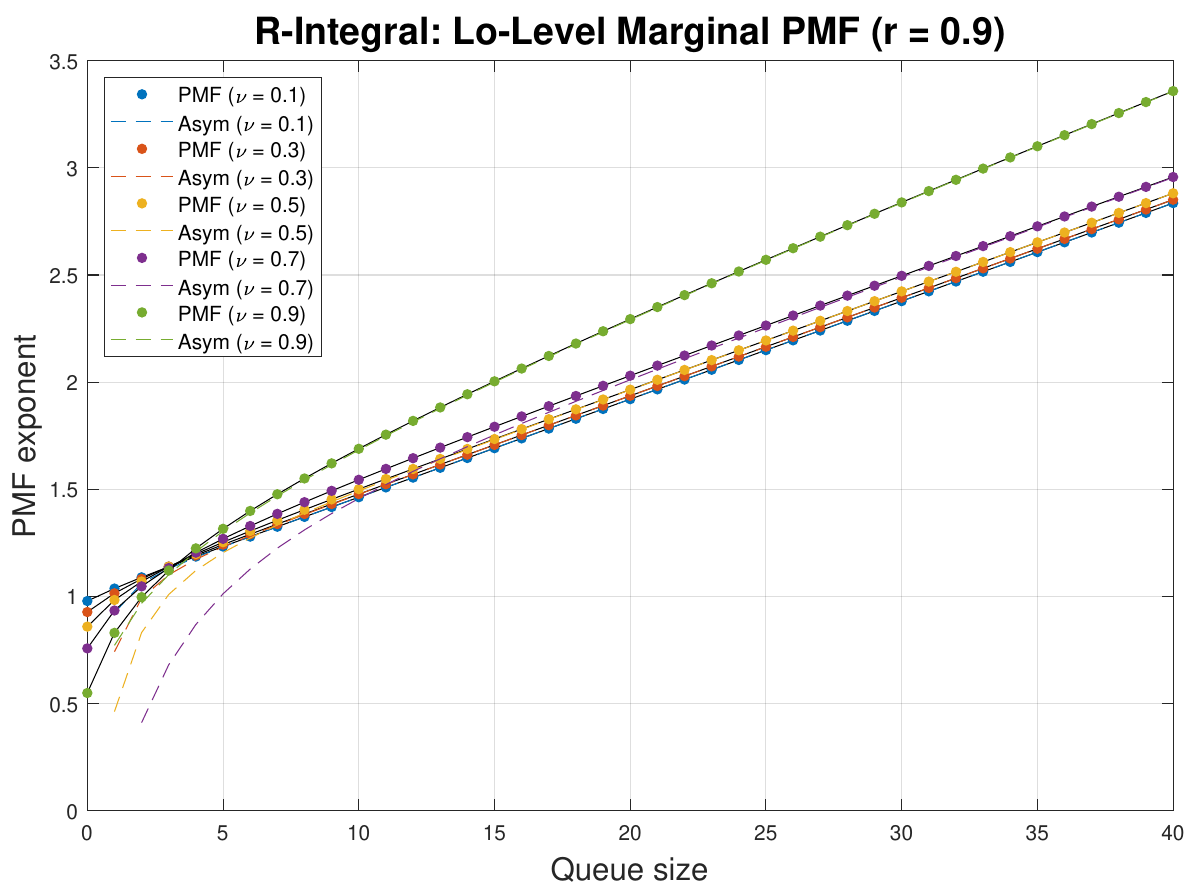}}
{\hphantom{x}\label{LoMargPMFMulti}}
{Queue-length PMF for the low-priority arrivals,
     plotted as the negative base-10 logarithm, for
     total traffic intensity $r = 0.9$ and a range of hifrac values ($\nu$).
     Asymptotic curves for the large queue-length limit are overlaid.}
\end{figure}

\begin{figure}
\FIGURE
{\includegraphics[width=\wscl\linewidth, height=\hscl\linewidth]{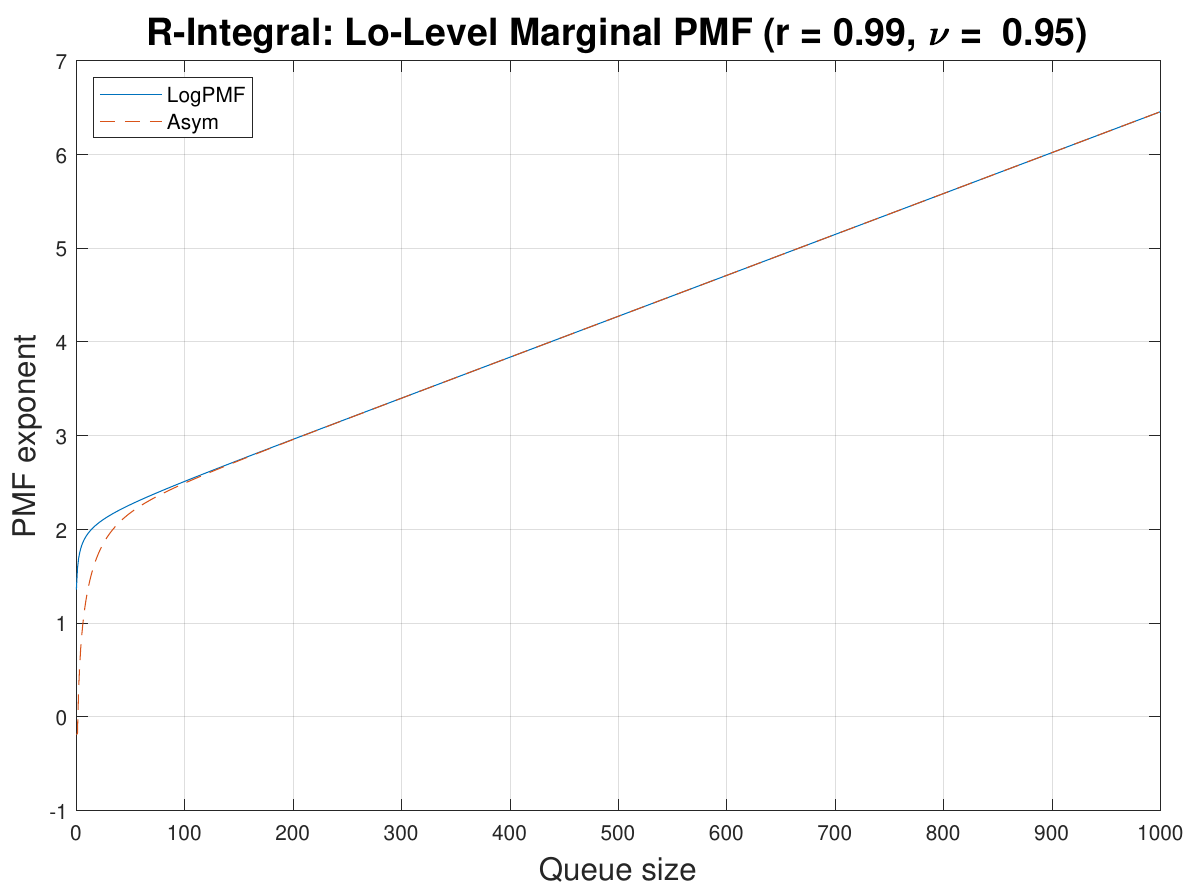}}
{\hphantom{x}\label{LoMargPMF}}
{Queue-length PMF for the low-priority arrivals,
     plotted as the negative base-10 logarithm, for
     total traffic intensity $r = 0.99$ and fraction of high-priority arrivals $\nu = 0.9$,
     with queue lengths extending far into the asymptotic region.
     It is compared with the exact asymptotic curve in the large queue-length limit.}
\end{figure}

\begin{figure}
\FIGURE
{\includegraphics[width=\wscl\linewidth, height=\hscl\linewidth]{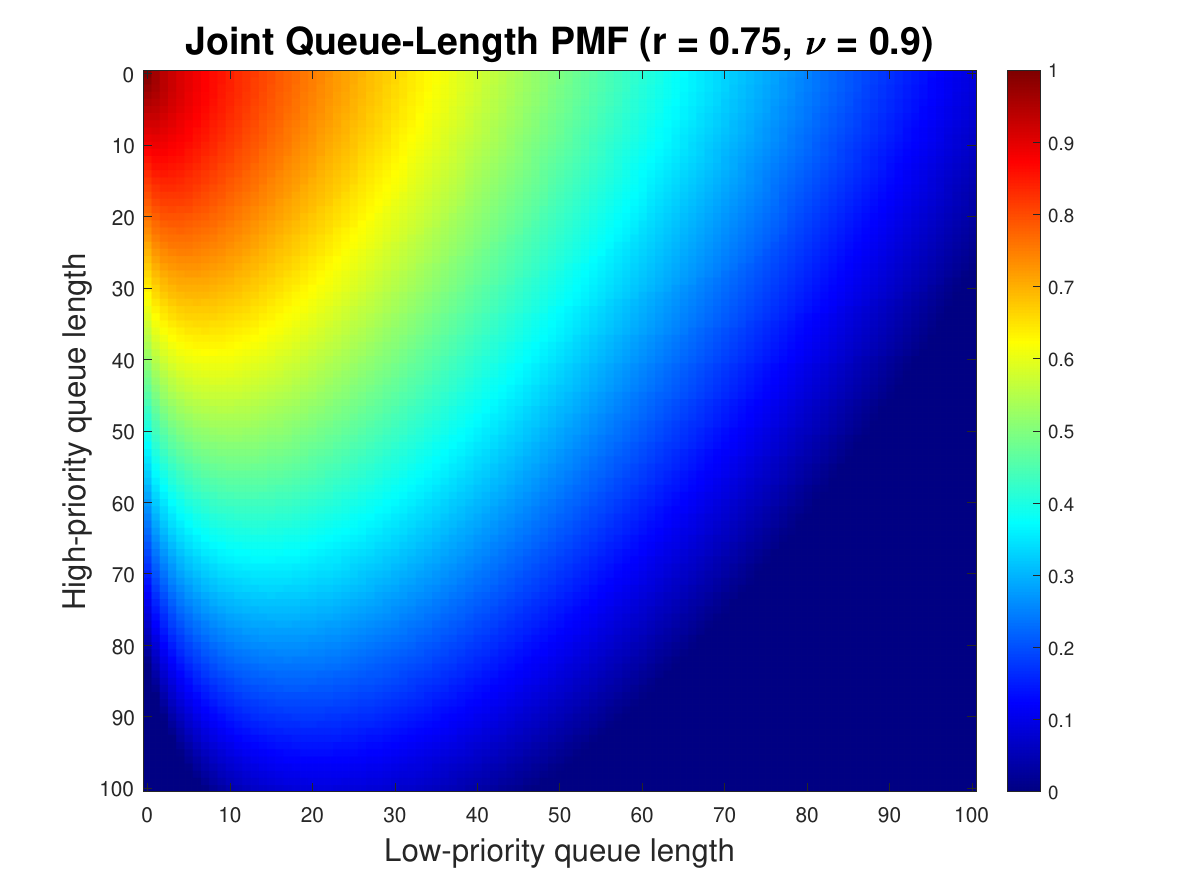}}
{\hphantom{x}\label{PMFMap}}
{2D map of the joint probability distribution of the queue lengths,
     with a logarithmic scaling applied,
     for total traffic intensity $r = 0.75$ and fraction of high-priority arrivals $\nu = 0.9$.}
\end{figure}

\section{Numerical Tests}
\begin{figure}
\FIGURE
{\includegraphics[width=\wsclx\linewidth, height=\hsclx\linewidth]{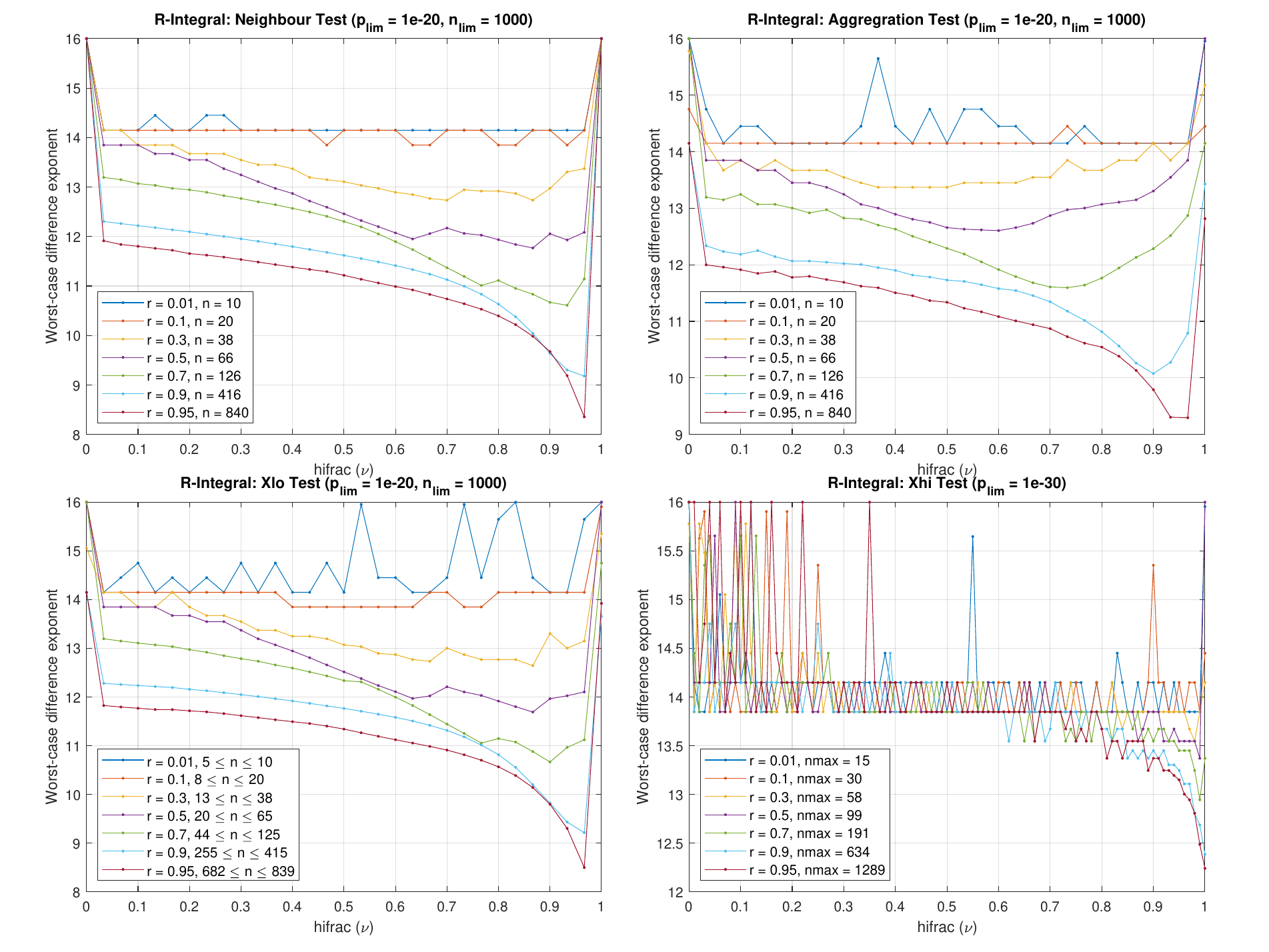}}
{\hphantom{x}\label{RTests}}
{Four tests of the joint probability distribution of the queue lengths as functions of the
     fraction of high-priority arrivals $\nu$, across a wide range of values
     for the total traffic intensity $r$, as displayed.
     The values on the vertical axes indicate the number of decimal places of agreement.}
\end{figure}

Various tests can be applied to quantify the numerical performance of the algorithm for
the computation of the joint PMF.

\subsection{Aggregation Test}
The aggregated queue-length distribution describes the total number of entities in the queue,
regardless of priority level.
This is equivalent to the queue-length distribution of the basic M/M/$c$ queueing model
with traffic intensity
\mbox{$r = r_{\text{lo}} + r_{\text{hi}}$},
which is known to be a simple geometric distribution.
Hence, the exact aggregate PMF is given by
\begin{equation}
f_{\text{agg}}^{(\text{ex})}(k) = (1-r)r^k \;,
\label{aggtest}
\end{equation}
for
\mbox{$k = 0,1,2,\ldots$}

One diagnostic test of the R-integral computational methodology is to check how well
the aggregate PMF constructed from the computed joint PMF reproduces the exact result.
This test is more convenient than similarly testing against the marginals as only a
finite summation is required.
Considering the joint PMF as a matrix whose rows and columns are labelled by its integer
arguments, values of the aggregate PMF are given by successive finite sums along the
anti-diagonals. Specifically,
in terms of the R-integrals, the aggregate PMF is expressed as
\begin{align}
\begin{aligned}
f_{\text{agg}}(k) &= \sum_{n=0}^k f(n,k-n) \\
&= \frac{1-r}{r}\sum_{n=0}^k(-r_{\text{lo}})^n(\Delta R_{n+1}^{k+2} - r_{\text{hi}}\Delta R_{n+1}^k) \;,
\end{aligned}
\end{align}
for
\mbox{$k = 0,1,2,\ldots$}.

We then consider the measure of performance (MOP)
\begin{equation}
\Xi_{\text{agg}} \equiv -\max_{k \geq 0}\left\{\log_{10}\left(|\ln(f_{\text{agg}}(k)) -
    \ln(f_{\text{agg}}^{(\text{ex})}(k))|\right)\right\} \;,
\end{equation}
where the maximum is taken over all values
\mbox{$0 \leq k \leq n_{\text{lim}}$}
such that
\mbox{$f_{\text{agg}}^{(\text{ex})}(k) > p_{\text{lim}} > 0$}.
Since we are working in double-precision arithmetic\footnote{All computation is performed in {\sc Matlab} R2020a,
which implements IEEE Standard 754 for double precision.},
all MOPs of this kind are capped
at a maximum allowed value of $16$.
The interpretation of $\Xi_{\text{agg}}$
(and similarly for all of the subsequent MOPs)
is that it indicates the number of decimal places of numerical agreement
in the worst case.

\subsection{Xhi-Test}
The exclusively-high distribution, defined by
\mbox{$f_{\text{xhi}}(m) \equiv f(0,m)$},
gives the probability of finding $m$ high-priority clients in the queue and no low-priority clients.
An exact expression for the exclusively-high probability is given by
\begin{equation}
f_{\text{xhi}}^{(\text{ex})}(m) = (1-r)(r_{\text{hi}}/z_2)^m  \;,
\label{xhitest}
\end{equation}
for
\mbox{$m = 0,1,2,\ldots$}
It is simple to calculate directly as the R-integral has only a simple pole when
\mbox{$n = 0$}.
One should note that $f_{\text{xhi}}(m)$ is not a proper PMF since
\mbox{$\sum_{m = 0}^\infty f_{\text{xhi}}(m) < 1$},
unless
\mbox{$\nu = 1$},
but can be turned into a conditional PMF by means of an overall scale factor.

In terms of the R-integrals, the exclusively-high PMF is expressed as
\begin{equation}
f_{\text{xhi}}(m) = \frac{1-r}{r}\left(\Delta R^{m+2}_{1} - z_1z_2\Delta R^{m}_{1}\right) \;,
\end{equation}
and we consider the MOP
\begin{equation}
\Xi_{\text{xhi}} \equiv -\max_{m \geq 0}\left\{\log_{10}\left(|\ln(f_{\text{xhi}}(m)) -
    \ln{(}f_{\text{xhi}}^{(\text{ex})}(m))|\right)\right\} \;,
\end{equation}
where the maximum is taken over all values
\mbox{$0 \leq m \leq n_{\text{lim}}$}
such that
\mbox{$f_{\text{xhi}}^{(\text{ex})}(m) > p_{\text{lim}} > 0$}.

\subsection{Xlo-Test}
Checking whether the computed joint PMF gives rise to the correct marginal distribution,
numerically, is not a convenient enterprise as it necessitates an infinite summation.
However, it is possible to devise an alternative test that checks the consistency of
the numerical low-priority marginal with the numerically computed joint PMF.
In the xlo-test, we relate the exclusively-low distribution $f_{\text{xlo}}(n)$
with the low priority marginal $f_{\text{lo}}(n)$.
To achieve this, we consider the PGF (\ref{gCohen}) recast into the form
\begin{equation}
g_m(p) = (1-r)\frac{1 - \lambda_1(p)}{\lambda_2(p) - r}{\cdot}\lambda_1^m(p) \;.
\end{equation}
Specialized to the case
\mbox{$m = 0$},
this may be expressed as
\begin{equation}
g_0(p) = (1-r)\left[1 + r_{\text{lo}}p{\cdot}\frac{1}{\lambda_2(p)-r}\right] \;.
\end{equation}
Since the PGF of the low-priority marginal is given by
\begin{equation}
g_{\text{lo}}(p) = \sum_{m=0}^\infty g_m(p) = \frac{1-r}{\lambda_2(p)-r} \;,
\end{equation}
we arrive at the result
\begin{equation}
g_0(p) = 1 - r + r_{\text{lo}}p{\cdot}g_{\text{lo}}(p) \;.
\end{equation}
There is a generalization of this result to non-zero values of $m$ that
relates $g_m(p)$ to $g_{\text{lo}}(p)$.
Its derivation is presented in the Appendix.
From the relationships
\begin{equation}
g_0(p) = \sum_{n=0}^\infty p^n f_{\text{xlo}}(n) \;, \quad
     g_{\text{lo}}(p) = \sum_{n=0}^\infty p^n f_{\text{lo}}(n) \;,
\end{equation}
we can equate powers to read off that
\begin{equation}
f_{\text{xlo}}(0) = 1 - r \;, \quad f_{\text{xlo}}(n) = r_{\text{lo}}{\cdot}f_{\text{lo}}(n-1) \;,
\end{equation}
for
\mbox{$n = 1,2,\ldots$},
or, equivalently,
\begin{equation}
f_{\text{xlo}}(n) = (1-r)\delta_{n0} + (1 - \delta_{n0})r_{\text{lo}}{\cdot}f_{\text{lo}}(n-1) \;,
\label{xlotest}
\end{equation}
for
\mbox{$n = 0,1,2,\ldots$},
where we can formally set
\mbox{$f_{\text{lo}}(-1) \equiv 0$}.
One should note that $f_{\text{xlo}}(n)$ is not a proper PMF since
\begin{equation}
\sum_{n= 0}^\infty f_{\text{xlo}}(n) = 1 - r_{\text{hi}} \;,
\end{equation}
but can be turned into a conditional PMF by means of an overall scale factor.

For the xlo-test, the LHS of (\ref{xlotest}) is taken to be given by (\ref{fxlo})
and is compared with the RHS of (\ref{xlotest}) where the marginal PMF
$f_{\text{lo}}(n)$
is expressed in terms of the R-integrals via (\refeq{flo}).
The relevant MOP is taken to be
\begin{equation}
\Xi_{\text{xlo}} \equiv -\max_{n > 0}\left\{\log_{10}\left(|\ln(f_{\text{xlo}}(n)) -
    \ln{(}r_{\text{lo}}f_{\text{lo}}(n-1))|\right)\right\} \;,
\end{equation}
where the maximum is taken over all values
\mbox{$0 < n \leq n_{\text{lim}}$}
such that
\mbox{$f_{\text{xlo}}(m) > p_{\text{lim}} > 0$}.

\subsection{Nearest-Neighbour Test}
A direct consequence of the recurrence relations for the R-integrals is that
the joint PMF at any given interior point
\mbox{$(n,m)$}
is a positively weighted sum of the joint PMF values at three of its four nearest neighbours:
\begin{equation}
f(n,m) = \frac{1}{1+r}\left[f(n,m+1) + r_{\text{lo}}f(n-1,m) + r_{\text{hi}}f(n,m-1)\right] \;,
\label{nntest}
\end{equation}
for all
\mbox{$m,n > 0$}.

In order to apply the neighbour test, we first compute the joint PMF
\mbox{$f(n,m)$}
on a 2D grid of points $(n,m)$ from (\ref{fjnt}).
Next, we use these values to compute the RHS of (\ref{nntest}), which we shall
denote $f_{\text{nn}}(n,m)$.
Then, we consider the MOP
\begin{equation}
\Xi_{\text{nn}} \equiv -\max_{m,n > 0}\left\{\log_{10}\left(|\ln(f(n,m)) -
     \ln(f_{\text{nn}}(n,m))|\right)\right\} \;,
\end{equation}
where the maximum is taken over all values
\mbox{$0 < m,n \leq n_{\text{lim}}$}
such that
\mbox{$f(n,m) > p_{\text{lim}} > 0$}.

\subsection{Quadratic Test}
In this test, we compare the results for the joint queue-length PMF computed from the R-integral
(denoted $f_{\text{ri}}(n,m)$) with that computed by the quadratic recurrence
(denoted $f_{\text{qr}}(n,m)$).
The MOP is taken to the be number of decimal places of agreement, as given by
\begin{equation}
\Xi_{\text{qr}} \equiv -\max_{m,n > 0}\left\{\log_{10}\left(|\ln(f_{\text{ri}}(n,m)) -
     \ln(f_{\text{qr}}(n,m))|\right)\right\} \;,
\end{equation}
where the maximum is taken over all values
\mbox{$0 < m,n \leq n_{\text{lim}}$}
such that
\mbox{$f(n,m) > p_{\text{lim}} > 0$}.

\subsection{Results}
Figure~\ref{RTests} presents the results of the numerical tests.
The MOP values relevant to the R-integral computations
are displayed on the vertical axis against the full range of
high-priority arrival fraction (hifrac) $\nu$ on the horizontal axis.
Individual curves are plotted for a discrete collection of traffic intensities,
spanning a wide range. Agreement always exceeds eight decimal places,
and is generally much higher.
The nearest-neighbour and xlo-tests check the internal consistency of the computations,
while the aggregation and xhi-tests check against exact analytical results.

The maximum queue occupancy to be examined was taken to be
\mbox{$n_{\text{lim}} = 1000$}.
PMF intervals examined included everything
down to a tail value of
\mbox{$p_{\text{lim}} = 10^{-20}$}
except in the xhi-test where
\mbox{$p_{\text{lim}} = 10^{-30}$}
was used.

Figure~\ref{QRecTest} presents the results of comparing the joint queue-length distribution
computed from the R-integral with that computed by the quadratic recurrence.
The close agreement observed implies a high level of accuracy for each method across
the complete range of parameters.
Worst case accuracy occurs when both the traffic intensity $r$ and hifrac $\nu$ approach unity.
In Table \ref{tab:acc}, we present results that investigate this region in more detail.
Values of $r$ close to unity have been reported to be problematic for the matrix-geometric
approach \citep{NP:Miller81,NP:Miller82}.
The table shows that both of the present methods behave well in this region.
The fourth and fifth columns indicate the smallest rectangular subset
\mbox{$[0,n_{\text{lo}}]\times [0,n_{\text{hi}}]$}
of
\mbox{$[0,n_{\text{lim}}]\times [0,n_{\text{lim}}]$}
that contains all grid points $(m,n)$ with probability greater than
\mbox{$p_{\text{lim}} = 10^{-20}$}.
A value of
\mbox{$n_{\text{lim}} = 1000$}
in one or both columns indicates that $p_{\text{lim}}$
was not attained in some direction. The last column is the minimum probability
that was achieved over all considered grid points whose probability values
exceed $p_{\text{lim}}$.
Computation time for the quadratic recurrence method is two orders of magnitude faster
than for the R-integral method.
Finally, Figure~\ref{LoMargQRecTest} repeats the quadratic test as described above,
but for the low-priority
marginal PMFs, with the distribution arising from the quadratic recurrence computed by
the algorithm of (\ref{SclRecur}).
The legend indicates the of range maximum queue lengths $n$ that had to be considered
across the full range of hifrac values $\nu$
in order to the attain the limiting probability level
\mbox{$p_{\text{lim}} = 10^{-20}$}
for the given traffic intensity $r$.
Agreement between the R-integral and quadratic recurrence approaches is observed
to exceed ten decimal places in the worst case..

\begin{table}
\renewcommand{\arraystretch}{1.3}
\TABLE
{Joint PMF Comparison\label{tab:acc}}
{\begin{tabular}{|c|c||c|c|c|c|}
\hline
$r$ & $\nu$ & $\Xi_{\text{qr}}$ & $n_{\text{hi}}$ &  $n_{\text{lo}}$ & $p_{\text{min}}$ \\
\hline\hline
$0.99$ & $0.95$  & $9.3279$  & $609$  & $1000$ & $1.0000\times 10^{-20}$ \\
       & $0.99$  & $8.1611$  & $1000$ & $1000$ & $1.0000\times 10^{-20}$ \\
       & $0.999$ & $6.6633$  & $1000$ & $1000$ & $1.0000\times 10^{-20}$ \\
       & $1.00$  & $11.7428$ & $1000$ & $0$    & $4.3171\times 10^{-7}$ \\
\hline
$0.999$ & $0.95$  & $9.4247$ & $685$  & $1000$ & $1.0000\times 10^{-20}$ \\
        & $0.99$  & $8.4169$ & $1000$ & $1000$ & $1.0017\times 10^{-20}$ \\
        & $0.999$ & $7.2251$ & $1000$ & $1000$ & $6.6926\times 10^{-18}$ \\
        & $1.00$  & $9.6972$ & $1000$ & $0$    & $3.6770\times 10^{-4}$ \\
\hline
$0.9999$ & $0.95$ & $9.4344$ & $657$  & $1000$ & $1.0000\times 10^{-20}$ \\
        & $0.99$  & $8.4361$ & $1000$ & $1000$ & $1.0000\times 10^{-20}$ \\
        & $0.999$ & $7.2455$ & $1000$ & $1000$ & $1.0540\times 10^{-18}$ \\
        & $1.00$  & $7.8504$ & $1000$ & $0$    & $9.0483\times 10^{-5}$ \\
\hline
\end{tabular}}
{}
\end{table}

\begin{figure}
\FIGURE
{\includegraphics[width=\wscl\linewidth, height=\hscl\linewidth]{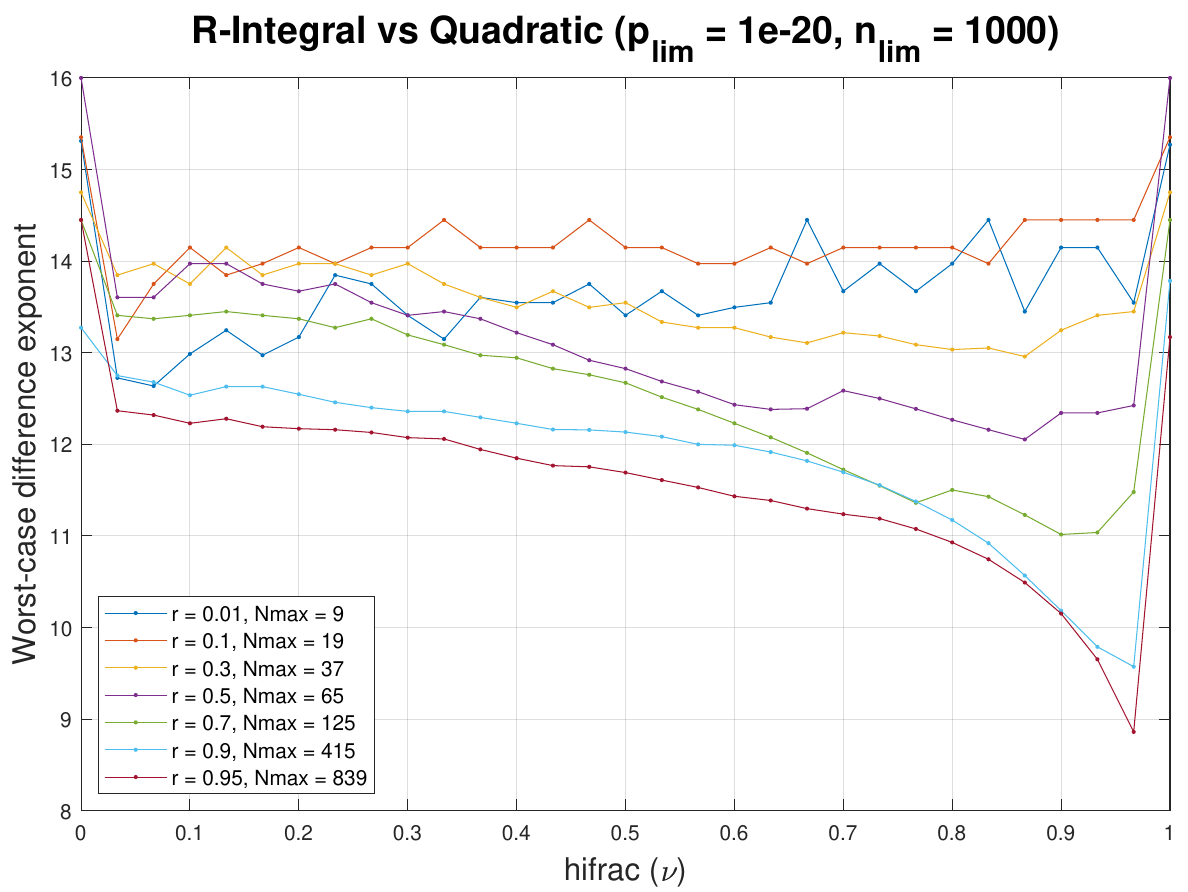}}
{\hphantom{x}\label{QRecTest}}
{Comparison of the joint probability distribution of the queue lengths as computed from
     the R-integral and from the quadratic recurrence.
     The MOP is plotted as a function of the
     fraction of high-priority arrivals $\nu$, across a wide range of values
     for the total traffic intensity $r$, as displayed.
     The values on the vertical axes indicate the number of decimal places of agreement.}
\end{figure}

\begin{figure}
\FIGURE
{\includegraphics[width=\wscl\linewidth, height=\hscl\linewidth]{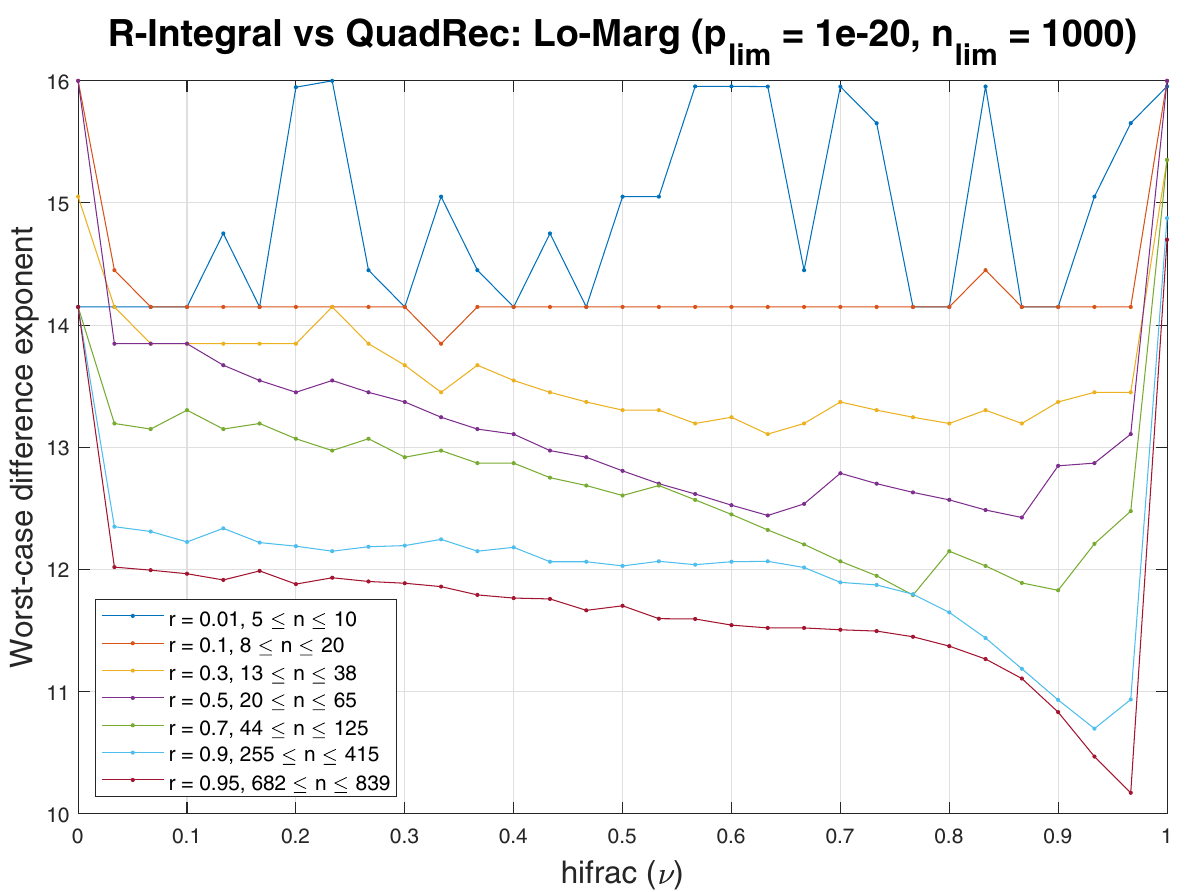}}
{\hphantom{x}\label{LoMargQRecTest}}
{Comparison of the low-priority marginal distribution of the queue lengths as computed from
     the R-integral and from the quadratic recurrence.
     The MOP is plotted as a function of the
     fraction of high-priority arrivals $\nu$, across a wide range of values
     for the total traffic intensity $r$, as displayed.
     The values on the vertical axes indicate the number of decimal places of agreement.}
\end{figure}

The exact results for the queue-lengths distributions derived here were also tested
against Monte-Carlo simulation. Excellent agreement was found across the entire parametric domain.
Details will be presented elsewhere.
We have also checked against the results in Table~3 of \citep{NP:Gail88} (where the service times are equal)
to find complete agreement. There, in the case of the present problem,
the quantity $P_{\text{Q}}$ is related to the no-wait probability
$P_{\text{NW}}$ given in (\ref{PNW}) by
\mbox{$P_{\text{Q}} = 1 - P_{\text{NW}}$},
and $p(0,0)$ is the probability that the system is empty, given by
\begin{equation}
\frac{1}{p(0,0)} = \frac{(rN)^N}{N!}\left[\frac{1}{1-r} + \Gamma_{\text{scl}}(rN,N)\right] \;,
\end{equation}
where
\begin{equation}
\Gamma_{\text{scl}}(x,\nu) \equiv \frac{\nu e^x}{x^\nu}\int_x^\infty dt\ t^{\nu-1}e^{-t}
\end{equation}
is the scaled upper incomplete gamma function as implemented in {\sc Matlab}.
In the present problem, neither of these quantities depend on the priority structure.
We relate the mean waiting times given in the table to the mean queue lengths via Little's law.

\section{Conclusions}
\label{Conclusions}
Simple methods for highly accurate computation of the joint and marginal distributions
for a non-preemptive two-level priority queue have been developed.
Explicit closed-form representations for the joint and marginal PMFs have also been derived,
something that has not been achieved previously.
Future work could entail extension of the present methods to unequal services rates
among the priority levels.

\begin{APPENDIX}{Convolutional Form}
In this appendix, we derive a relationship between joint PMF
\mbox{$f(n,m)$}
and the low-priority marginal
\mbox{$f_{\text{lo}}(n)$}.
One may observe the general structure
\begin{equation}
g_m(p) = A^{(m)}(p) + B^{(m)}(p){\cdot}g_{\text{lo}}(p) \;,
\end{equation}
for polynomials
\begin{equation}
A^{(m)}(p) = \sum_{n=0}^m A^{(m)}_n p^n \;, \quad B^{(m)}(p) = \sum_{n=0}^m B^{(m)}_n p^n \;.
\end{equation}
Thus, we obtain the convolutional form for the joint PMF:
\begin{equation}
f(n,m) = A^{(m)}_n + \sum_{k=0}^n B^{(m)}_k{\cdot} f_{\text{lo}}(n-k) \;,
\end{equation}
which generalizes the relationship between
\mbox{$f_{\text{xlo}}$}
and
\mbox{$f_{\text{lo}}$}
given in (\ref{xlotest})
to non-zero values of $m$.
One may also note the special case
\begin{equation}
f_{\text{xhi}}(m) \equiv f(0,m) = A^{(m)}_0 + B^{(m)}_0 f_{\text{lo}}(0) \;.
\end{equation}
In what follows, we derive explicit expressions for
$A^{(m)}(p)$ and $B^{(m)}(p)$.

Since
\mbox{$\lambda_1(p)$}
satisfies a quadratic equation, we have that
\mbox{$\lambda_1^m(p) = \alpha_m(p)\lambda_1(p) + \beta_m(p)$}
for some polynomials
\mbox{$\alpha_m(p)$},
\mbox{$\beta_m(p)$}.
The fact that $\alpha_m,\beta_m$ are polynomials follows from examining
\mbox{$\lambda_1^2$}.
Let us now recall that the Chebyshev polynomials of the first and second kind,
$T_n(x)$ and $U_n(x)$ respectively,
may be expressed as
\begin{align}
\begin{aligned}
T_n(x)                 &= \half\left[(x + \sqrt{x^2 - 1})^n + (x - \sqrt{x^2 - 1})^n\right] \;, \\
\sqrt{x^2-1}U_{n-1}(x) &= \half\left[(x + \sqrt{x^2 - 1})^n - (x - \sqrt{x^2 - 1})^n\right] \;,
\end{aligned}
\end{align}
which implies the identity
\begin{equation}
(x - \sqrt{x^2 - 1})^n = T_n(x) - xU_{n-1}(x) + (x - \sqrt{x^2 - 1}){\cdot}U_{n-1}(x) \;.
\end{equation}
Since we have
\begin{equation}
\lambda_1(p) = b(p) - \sqrt{b^2(p) - r_1} \;, \quad
     b(p) \equiv (1 + r - r_2 p)/2 \;,
\end{equation}
it follows that
\begin{equation}
\frac{\lambda_1^m(p)}{r_1^{m/2}} = T_m(x(p)) - x(p){\cdot}U_{m-1}(x(p)) +
     \frac{\lambda_1(p)}{\sqrt{r_1}}{\cdot}U_{m-1}(x(p)) \;,
\end{equation}
with
\mbox{$x(p) \equiv b(p)/\sqrt{r_1}$}.
Consequently,
\begin{align}
\begin{aligned}
\alpha_m(p) &= r_1^{(m-1)/2}U_{m-1}(x(p)) \;, \\
\beta_m(p)  &= r_1^{m/2}\left[T_m(x(p)) - x(p){\cdot}U_{m-1}(x(p))\right] \;,
\label{alphabeta}
\end{aligned}
\end{align}
for
\mbox{$x = 0,1,\ldots$},
with
\mbox{$U_{-1}(x) \equiv 0$}.
The first few $\alpha$-coefficients are given by
\begin{equation}
\alpha_0(p) = 0 \;, \quad \alpha_1(p) = 1 \;, \quad  \alpha_2(p) = 1 + r - r_2 p \;.
\end{equation}
The first few $\beta$-coefficients are given by
\begin{equation}
\beta_0(p) = 1 \;, \quad \beta_1(p) = 0 \;, \quad  \beta_2(p) = -r_1 \;.
\end{equation}
Using identities satisfied by the Chebyshev polynomials, (\ref*{alphabeta}) can be simpified as
\begin{equation}
\alpha_m(p) = r_1^{(m-1)/2}U_{m-1}(x(p)) \;, \quad
     \beta_m(p)  = -r_1^{m/2}U_{m-2}(x(p)) \;,
\label{alphabeta1}
\end{equation}
which implies the relationship
\mbox{$\beta_m(p) = -r_1\alpha_{m-1}(p)$},
where we formally set
\mbox{$\alpha_{-1}(p) \equiv 1$}.

Noting that
\begin{align}
\begin{aligned}
(1 - \lambda_1)\lambda_1^m &= (\alpha_m - \alpha_{m+1})\lambda_1 + (\beta_m - \beta_{m+1}) \\
     &= -(\alpha_m - \alpha_{m+1})(\lambda_2-r) + (1-r_2p)(\alpha_m - \alpha_{m+1}) + (\beta_m - \beta_{m+1}) \;,
\end{aligned}
\end{align}
followed by substitution into the representation
\begin{equation}
g_m(p) = \frac{1 - r}{\lambda_2-r}{\cdot}(1- \lambda_1)\lambda_1^m \;,
\end{equation}
yields
\begin{equation}
g_m(p) = -(1-r)(\alpha_m - \alpha_{m+1}) + g_{\text{lo}}(p)\left[
     (1-r_2p)(\alpha_m - \alpha_{m+1}) + (\beta_m - \beta_{m+1})\right] \;,
\end{equation}
from which we can read off
\begin{align}
\begin{aligned}
A^{(m)}(p) &= (1-r)\left[\alpha_{m+1}(p) - \alpha_m(p)\right] \;, \\
B^{(m)}(p) &= (1 - r_2p)\left[\alpha_m(p) - \alpha_{m+1}(p)\right] +
     \left[\beta_m(p) - \beta_{m+1}(p)\right] \;.
\end{aligned}
\end{align}
In terms of Chebyshev polynomials, this becomes
\begin{align}
\begin{aligned}
A^{(m)}(p) &= -(1-r)r_1^{(m-1)/2}\left[U_{m-1}(x(p)) - \sqrt{r_1}U_{m}(x(p))\right] \;, \\
B^{(m)}(p) &= -\frac{1-r_2p}{1-r}{\cdot}A^{(m)}(p) + \frac{r_1}{1-r}{\cdot}A^{(m-1)}(p)\;,
\end{aligned}
\end{align}
for
\mbox{$m = 0,1,\ldots$},
where we set
\mbox{$U_{-1}(x) \equiv 0$},
\mbox{$U_{-2}(x) \equiv -1$}.
As a sanity check, it is straightforward to confirm that
\mbox{$A^{(m)}(1) + B^{(m)}(1) = (1-r_1)r_1^m$}.
Another check is given by
\begin{align}
\begin{aligned}
\sum_{m=0}^\infty A^{(m)}(p)&= \left\{
\begin{array}{cll}
(1 - r)/(1 - r_1) & \quad\text{for}\quad & p = 1 \\
0                 & \quad\text{for}\quad & p \neq 1
\end{array}
\right. \;, \\
\sum_{m=0}^\infty B^{(m)}(p) &= \left\{
\begin{array}{cll}
(r - r_1)/(1 - r_1) & \quad\text{for}\quad & p = 1 \\
1                   & \quad\text{for}\quad & p \neq 1
\end{array}
\right. \;. \\
\end{aligned}
\end{align}
\end{APPENDIX}


%
\section*{Acknowledgments.}
The authors gratefully acknowledge useful discussions with Dr.~Stephen Bocquet.


\bibliographystyle{informs2014} 
\bibliography{NPPriority} 



\end{document}